\documentclass[11pt]{amsart}
\usepackage{times}
\usepackage[mathscr]{eucal}
\usepackage{epsfig}
\usepackage{amssymb}

\newcommand{\hcase}[2]%
{\makebox[0pt]%
{\raisebox{-1pt}[0pt][0pt]{#1{#2}}}}%

\newcommand{\hbicase}[2]%
{\makebox[0pt]%
{\raisebox{-2.5pt}[0pt][0pt]{#1{#2}}}}%

\newcommand{\EAR}[1]%
{\begin{picture}(#1,0)%
\put(0,0){\vector(1,0){#1}}%
\end{picture}}%

\newcommand{\EEQL}[1]%
{\begin{picture}(#1,0)%
\truex{200}%
\put(0,\value{x}){\line(1,0){#1}}%
\put(0,0){\line(1,0){#1}}%
\end{picture}}%

\newcommand{\ear}%
{\hspace{\SOURCE\unitlength}%
\hcase{\EAR}{\ARROWLENGTH}}%

\newcommand{\eeql}%
{\hspace{\SOURCE\unitlength}%
\hbicase{\EEQL}{\ARROWLENGTH}}%

\newcommand{\vcase}[2]{#1{#2}}%

\newcommand{\vbicase}[2]{\makebox[0pt]{{#1{#2}}}}%

\newcommand{\SAR}[1]%
{\begin{picture}(0,0)%
\put(0,0){\makebox(0,0)%
{\begin{picture}(0,#1)%
\put(0,#1){\vector(0,-1){#1}}%
\end{picture}}}\end{picture}}%

\newcommand{\SEQL}[1]%
{\begin{picture}(0,0)%
\truex{100}%
\put(0,0){\makebox(0,0)%
{\begin{picture}(0,#1)\put(-\value{x},#1){\line(0,-1){#1}}%
\put(\value{x},#1){\line(0,-1){#1}}%
\end{picture}}}\end{picture}}%

\newcommand{\sarv}[1]{\vcase{\SAR}{#100}}%

\newcommand{\sar}{\sarv{50}}%

\newcommand{\seqlv}[1]{\vbicase{\SEQL}{#100}}%

\newcommand{\seql}{\seqlv{50}}%

\newcount\SCALE%

\newcount\NUMBER%

\newcount\LINE%

\newcount\COLUMN%

\newcount\WIDTH%

\newcount\SOURCE%

\newcount\ARROW%

\newcount\TARGET%

\newcount\ARROWLENGTH%

\newcount\NUMBEROFDOTS%

\newcounter{x}%

\newcounter{y}%

\newcounter{z}%

\newcounter{horizontal}%

\newcounter{vertical}%

\newskip\itemlength%

\newskip\firstitem%

\newskip\seconditem%

\newcommand{\printarrow}{}%

\newcommand{\truex}[1]{%
\NUMBER=#1%
\multiply\NUMBER by 100%
\divide\NUMBER by \SCALE%
\setcounter{x}{\NUMBER}}%

\newcommand{\changecounters}[1]{%
\SOURCE=\ARROW%
\ARROW=\TARGET%
\settowidth{\itemlength}{#1}%
\ifdim \itemlength > 2800\unitlength%
\addtolength{\itemlength}{-2800\unitlength}%
\TARGET=\itemlength%
\divide\TARGET by 1310%
\multiply\TARGET by 100%
\divide\TARGET by \SCALE%
\else%
\TARGET=0%
\fi%
\ARROWLENGTH=5000%
\advance\ARROWLENGTH by -\SOURCE%
\advance\ARROWLENGTH by -\TARGET%
\advance\SOURCE by -\TARGET}%

\newcommand{\initialize}[1]{%
\LINE=0%
\COLUMN=0%
\WIDTH=0%
\ARROW=0%
\TARGET=0%
\changecounters{#1}%
\renewcommand{\printarrow}{#1}%
\begin{center}%
\vspace{10pt}%
\begin{picture}(0,0)}%

\newcommand{\DIAGV}[2]{%
\SCALE=#1%
\setlength{\unitlength}{655sp}%
\multiply\unitlength by \SCALE%
\divide\unitlength by 100%
\initialize{\mbox{$#2$}}}%

\newcommand{\n}[1]{%
\changecounters{\mbox{$#1$}}%
\put(\COLUMN,\LINE){\makebox(0,0){\printarrow}}%
\thinlines%
\renewcommand{\printarrow}{\mbox{$#1$}}%
\advance\COLUMN by 4000}%

\newcommand{\nn}[1]{%
\put(\COLUMN,\LINE){\makebox(0,0){\printarrow}}%
\thinlines%
\ifnum \WIDTH < \COLUMN%
\WIDTH=\COLUMN%
\else%
\fi%
\advance\LINE by -4000%
\COLUMN=0%
\ARROW=0%
\TARGET=0%
\changecounters{\mbox{$#1$}}%
\renewcommand{\printarrow}{\mbox{$#1$}}}%

\newcommand{\conclude}{%
\put(\COLUMN,\LINE){\makebox(0,0){\printarrow}}%
\thinlines%
\ifnum \WIDTH < \COLUMN%
\WIDTH=\COLUMN%
\else%
\fi%
\setcounter{horizontal}{\WIDTH}%
\setcounter{vertical}{-\LINE}%
\end{picture}}%

\newcommand{\diag}{%
\conclude%
\raisebox{0pt}[0pt][\value{vertical}\unitlength]{}%
\hspace*{\value{horizontal}\unitlength}%
\vspace{10pt}%
\end{center}%
\setlength{\unitlength}{1pt}}%

\newcommand{\diagv}[3]{%
\conclude%
\NUMBER=#1%
\rule{0pt}{\NUMBER pt}%
\hspace*{-#2pt}%
\raisebox{0pt}[0pt][\value{vertical}\unitlength]{}%
\hspace*{\value{horizontal}\unitlength}
\NUMBER=#3%
\advance\NUMBER by 10%
\vspace*{\NUMBER pt}%
\end{center}%
\setlength{\unitlength}{1pt}}%

\newtheorem{theorem}{Theorem}[section]
\newtheorem{corollary}[theorem]{Corollary}
\newtheorem{lemma}[theorem]{Lemma}
\newtheorem{proposition}[theorem]{Proposition}

\theoremstyle{definition}
\newtheorem{definition}[theorem]{Definition}
\theoremstyle{remark}

\numberwithin{equation}{section}

\def\Hom{\mbox{\rm Hom}}
\def\Ext{\mbox{\rm Ext}}

\def\Qco{\mathfrak{Qco}}
\def\Z{\mathbb{Z}}
\def\Lt{{\mathcal{L}}}
\def\Lto{{\mathcal{L}}^{\perp}}
\def\Ltoo{{ ^{\perp}{\mathcal{L}}}}
\def\F{\mathcal{F}}

\def\E{{\mathscr E}}
\def\O{{\mathcal{O}}}
\def\R{{\mathfrak R}}
\def\C{{\mathcal{C}}}
\def\J{{\mathcal{J}}}

\begin{document}

\title{Gorenstein Categories and Tate Cohomology on Projective Schemes}

\address{Department of Mathematics, University of Kentucky, Lexington, Kentucky
40506-0027, U.S.A.} \email{enochs@ms.uky.edu}

\author{E. Enochs, S. Estrada and J.R.Garc\'{\i}a Rozas}
\address{
Departamento de Matem\'atica Aplicada, Universidad de Murcia,Campus
del Espinardo, Espinardo (Murcia) 30100, Spain}
\email{sestrada@um.es}%
\address{
Departamento de \'Algebra y A. Matem\'atico, Universidad de
Almer\'{\i}a, Almer\'{\i}a 04071, Spain}\email{jrgrozas@ual.es}

\thanks{The authors are partially supported by the DGI MTM2005-03227}

\date{}

\begin{abstract}
We study Gorenstein categories. We show that such a category has
Tate cohomological functors and Avramov-Martsinkovsky exact
sequences connecting the Gorenstein relative, the absolute and the
Tate cohomological functors. We show that such a category has what
Hovey calls an injective model structure and also a projective model
structure in case the category has enough projectives.\\
As examples we show that if $X$ is a locally Gorenstein projective
scheme then the category $\Qco(X)$ of quasi-coherent sheaves on $X$
is such a category and so has these features.
\end{abstract}
\maketitle

\section{Introduction}
Tate homology and cohomology over the group ring $\Z G$ (with $G$ a
finite group) began with Tate's observation that the $\Z G$-module
$\Z$ with the trivial action admits a complete projective
resolution. Apparently motivated by Tate's work, Auslander showed in
\cite{AB} that if $A$ is what Bass in \cite{bass} calls a Gorenstein
local ring, the finitely generated maximal Cohen-Macaulay modules
can be characterized as those which admit a complete projective
resolution of finitely generated projective modules. Auslander calls
these modules the modules of $G$-dimension $0$ and goes on to define
the $G$-dimension of any finitely generated module.

In \cite{EnOv} an easy modification of one of Auslander's
characterizations of the finitely generated modules of $G$-dimension
$0$ was given and so allows one to extend the definition to any
module (finitely generated or not). Since this modified definition
dualizes it seems appropriate to use the terms Gorenstein projective
(corresponding to modules of $G$-dimension $0$) and Gorenstein
injective. Then there is a natural way to define Gorenstein
projective and injective dimension of any module. If a module $M$
has finite Gorenstein projective dimension $n$, then the $n$-th
syzygy of $M$ has a complete projective resolution. This complex is
a homotopy invariant of $M$ and so can be used to define Tate
homological functors $\widehat{\Ext}_R^n(M,N)$ and $\widehat{{\rm
Tor}}_n^R(M,N)$ for any $n\in \Z$. If on the other hand $N$ has
finite Gorenstein injective dimension a similar procedure can be
used to define analogous functors. A. Iacob in \cite{iacob2} showed
that if both conditions hold then the two procedures give us the
same functors, i.e., that we have balance in this situation.

In categories of sheaves there are usually not enough projectives
but there are enough injectives. So we use the second approach to
define Gorenstein injectives on the category of quasi-coherent
sheaves on certain projective schemes. More precisely, we show that
if such a scheme $X\subseteq {\bf P}^n(A)$ (where $A$ is
commutative noetherian) is what we call a locally Gorenstein scheme
then all objects of $\Qco(X)$ have finite Gorenstein injective
dimension and that there is a universal bound of these dimensions.
This allows us to define Tate cohomology in this situation, to get
Avramov-Martsinkovsky sequences and to impose a model structure on
$\Qco(X)$.

The example $\Qco(X)$ mentioned above is our motivation for defining
Gorenstein categories. These categories will be Grothendieck
categories with properties much like those of categories of modules
over Gorenstein rings. But we would like our definitions to be such
that nice categories of sheaves will be Gorenstein. Since categories
of sheaves rarely have enough projectives we need a definition which
does not involve projective objects. But such categories do have
enough injectives and so the functors $\Ext$ are defined. And so
projective dimensions of objects can be defined in terms of the
vanishing of the $\Ext$ functors. So we define a Gorenstein category
in terms of the global finitistic projective and injective
dimensions of the category. After defining Gorenstein categories and
proving some basic results about them, we consider examples. We show
that if $X\subset {\bf P}^n(A)$ (where $A$ is commutative
noetherian) is a projective scheme $\Qco(X)$ will be a Gorenstein
category when $X$ is a locally Gorenstein scheme.

The authors give a sincere thanks to Mark Hovey for his interest in
our work and for motivating us by his work on cotorsion pairs and
model structures.

For all unexplained terminology see \cite{EdO}.

\section{ Gorenstein Categories}

Our object now is to define Gorenstein categories and then exhibit
some of their properties. In the section ${\mathcal{A}}$ will always
be a Grothendieck category with a specified generator $X$. The
symbols $Y,Z$ etc. will denote objects of ${\mathcal{A}}$. We will
use the generator $X$ to assign a cardinal number to every object
$Y$. This cardinal will be $|\Hom(X,Y)|$. We will now give several
lemmas with the object of showing that we can do what is called
set-theoretic homological algebra in a Grothendieck category.

We refer the readers to \cite{BoSten} for the definition and basic
properties of a Grothendieck category.
\begin{lemma}
If $I$ is a set and if $X^{(I)}\rightarrow Y$ is an epimorphism,
then there is a subset $J\subset I$ with $|J|\leq |Y|$ and such that
the restriction $X^{(J)}\rightarrow Y$ is also an epimorphism.
\end{lemma}
\begin{proof} For such an $X^{(I)}\rightarrow Y$ we see that for
each $i\in I$ we have an associated morphism $X\rightarrow Y$. Let
$J\subset I$ be such that the morphisms $X\rightarrow Y$
corresponding to the $j\in J$ give all these morphisms and such that
if $j\neq j'$ ($j,j'\in J$) then $j$ and $j'$ correspond to
different morphisms.  Then $|J|\leq |\Hom (X,Y)|=|Y|$. Also there is
a natural factorization $X^{(I)}\rightarrow X^{(J)}\rightarrow Y$.
So $X^{(J)}\rightarrow Y$ is also an epimorphism.
\end{proof}

\begin{lemma}
For any $Y$ there is an epimorphism $X^{(|Y|)}\rightarrow Y$.
\end{lemma}
\begin{proof} Immediate. \end{proof}
\begin{corollary}\label{epicard}
Given any object $Y$ of ${\mathcal{A}}$ there is a cardinal $\kappa$
such that if $U\rightarrow Y$ is an epimorphism then there is a
subobject $U'\subset U$ such that $U'\rightarrow Y$ is an
epimorphism and such that $|U'|\leq\kappa$.
\end{corollary}
\begin{proof} Any such $U$ is a quotient of $X^{(I)}$ for some set
$I$. Then $X^{(I)}\rightarrow U\rightarrow Y$ is an epimorphism. But
from the above we see that there is a subset $J\subset I$ with
$|J|\leq |Y|$ such that $X^{(J)}\rightarrow U\rightarrow Y$ is an
epimorphism. So let $U'$ be the image of $X^{(J)}$ in $U$. So we see
that it is easy to get a $\kappa$ that bounds all $U'$ that we get
in this manner. \end{proof}

\begin{corollary}\label{setgen}
For every cardinal $\kappa$ there is a set of representatives of
objects $Y$ with $|Y|\leq \kappa$.
\end{corollary}
\begin{proof} By having such a set we mean that we have a set of
$Y$ with $|Y|\leq \kappa$ such that every $Z$ with $|Z|\leq \kappa$
is isomorphic to some $Y$ in our set. Now let $|Y|\leq \kappa$. Then
from the above $Y$ is the quotient of the coproduct $X^{(\kappa)}$.
But we can clearly form a set of representatives of quotients of
$X^{(\kappa)}$ since ${\mathcal{A}}$ is well-powered (see
\cite[Proposition 10.6.3]{libcat}).
\end{proof}

\begin{lemma}\label{bob}
If $Z\subset Y$ then $|Z|\leq |Y|$.
\end{lemma}
\begin{proof} Immediate. \end{proof}

\begin{lemma}\label{cociente}
Given a cardinal $\kappa$ there exists a cardinal $\lambda$ such
that if $|Y|\leq \kappa$ and if $Z\subset Y$ then $|Y/Z|\leq
\lambda$.
\end{lemma}
\begin{proof} By the Corollary above there is an epimorphism
$X^{(|Y|)}\rightarrow Y$ and so an epimorphism
$X^{(\kappa)}\rightarrow Y$. Hence $Y$ is isomorphic to a quotient
of $X^{(\kappa)}$. Since there is a set of representatives of
subobjects $S$ of $X^{(\kappa)}$, we can take $\lambda$ to be the
sup of all $|X^{(\kappa)}/S|$ taken over all $S\subset
X^{(\kappa)}$.  \end{proof}

\begin{lemma}\label{disum}
For any object $Y$ and set $I$ we have $|Y^{(I)}|\leq
|Y^I|=|Y|^{|I|}$.
\end{lemma}
\begin{proof} The equality follows from the equality
$\Hom(X,Y^I)=\Hom(X,Y)^I$. Since we are in a Grothendieck category $
Y^{(I)}\to Y^I$ is a monomorphism, so we get the inequality.
\end{proof}

\begin{lemma}
For any objects $Y$ and $Z$ we have $|\Hom(Y,Z)|\leq |Z|^{|Y|}$.
\end{lemma}
\begin{proof} We have an epimorphism $X^{(|Y|)}\rightarrow Y$ so
$\Hom(Y,Z)\subset \Hom(X^{(|Y|)},Z)=\Hom(Y,Z)^{|Y|}$ But
$|\Hom(Y,Z)^{|Y|}|=|Z|^{|Y|}$. \end{proof}

\begin{lemma}
If $\gamma$ is an ordinal and if $(\kappa_{\alpha})_{\alpha
<\gamma}$ is a family of cardinal numbers, then there is a cardinal
number $\lambda$ such that if $(Y_{\alpha})_{\alpha < \gamma}$ is a
family of objects with $Y_{\alpha}\subset Y_{\alpha'}$ whenever
$\alpha \leq \alpha' , \gamma$ and such that $|Y_{\alpha}|\leq
\kappa_{\alpha}$ for each $\alpha <\gamma$ then $|\cup
Y_{\alpha}|\leq \lambda$.
\end{lemma}
\begin{proof} We have an epimorphism $X^{(|Y_{\alpha}|)}\rightarrow
Y_{\alpha}$ for each $\alpha$. So since we are in a Grothendieck
category we have an epimorphism $X^{(\sum |Y_{\alpha}|)}\rightarrow
\cup Y_{\alpha}$ and so we have an epimorphism $X^{(\sum
\kappa_{\alpha})}\rightarrow \cup Y_{\alpha}$.
 So $\cup Y_{\alpha}$ is a quotient object of $X^{(\sum \kappa_{\alpha})}$.
So now we appeal to Lemma \ref{cociente} and get our $\lambda$.
\end{proof}

\begin{lemma}
Given a cardinal $\kappa$ there is a cardinal $\lambda$ such that if
$|Y|\leq \kappa$ then $|E(Y)|\leq \lambda$ (here $E(Y)$ is an
injective envelope of $Y$).
\end{lemma}
\begin{proof} We only need argue that for a given $\kappa$ there is
a $\lambda$ such that if $|Y|\leq \kappa$ then $Y\subset E$ for an
injective object $E$ where $|E|\leq \lambda$. To show this we use a
categorical version of Baer's original proof that every module can
be embedded in an injective module. Since $X$ is a generator Baer's
criterion says an object $E$ is injective if and only if it is
injective for $X$, i.e. if for any subobject $S\subset X$ and any
morphism $S\rightarrow E$ there is an extension $X\rightarrow E$.
Given the object $Y=Y_0$ Baer first constructs $Y_1$ with
$Y_0\subset Y_1$ and such that for any $S\subset X$ and any
$S\rightarrow Y_0$ there is an extension $S\rightarrow Y_1$. By his
construction $Y_1$ is the quotient of the coproduct of $Y_0$ and
copies of $X$ where this is a copy of $X$ for each $S\rightarrow
Y_0$ (with $S\subset X$ arbitrary). The quotient identifies each
such $S\subset X$ for a given $S\rightarrow Y_0$ with its image in
$Y_0$. By Lemmas \ref{cociente} and \ref{disum} we see that if
$\kappa =\kappa_0$ is a cardinal we can find a cardinal $\kappa_1$
such that if $|Y|=|Y_0|\leq \kappa$ then the $Y_1$ as constructed
above is such that $|Y_1|\leq \kappa_1$. Then for any ordinal
$\beta$ Baer constructs a continuous chain $(Y_{\alpha})_{\alpha
\leq \beta }$ of objects (so $Y_{\alpha} \subset Y_{\alpha'}$ if
$\alpha\leq \alpha' \leq \beta$ and if $\gamma \leq \beta$ is a
limit ordinal then $Y_{\gamma}=\cup Y_{\alpha}$ (${\alpha <
\gamma}$) where $Y_0$ is a given $Y$. And then if $\alpha +1\leq
\beta $ then we get $Y_{\alpha +1}$ for $Y_{\alpha}$ just as $Y_1$
is constructed from $Y_0$ as above and where $Y_{\gamma}=\cup
Y_{\alpha}$ ${\alpha < \gamma}$. Then for a given $\kappa$ we see
that we can find a family $(\kappa_{\alpha})_{\alpha \leq \beta} $
of cardinal numbers such that $\kappa_0 =\kappa$ and such that if
$(Y_{\alpha})_{\alpha \leq \beta}$ is constructed as above where
$|Y|=|Y_0|\leq \kappa$ then $|Y_{\alpha}|\leq \kappa_{\alpha}
$ for each $\alpha \leq \beta$.\\
 We now note that if $S\subset X$ and if $S\rightarrow Y_{\beta}$ is such
that there is a factorization $S\rightarrow Y_{\alpha}\rightarrow
Y_{\beta}$ for some $\alpha <\beta $ then $S\rightarrow Y_{\alpha}$
has an extension $X\rightarrow Y_{\alpha +1}$ so giving the
extension $X\rightarrow Y_{\alpha +1}\rightarrow Y_{\beta}$ of the
original $S\rightarrow Y_{\beta}$. So we must choose $\beta$ such
that every such $S\rightarrow Y_{\beta}$ has such a factorization.
So we want to argue that we can find a $\beta$ independent of $Y$ so
that the corresponding $Y_{\beta}$ will always be injective. But
this again just uses Baer's original idea and appeals to the fact
that every object is small relative to the class of monomorphisms
(see [22, pg. 32] for the terminology and [1, Proposition 2.2] for
an argument). The object in question would be the coproduct of a
representative set of subobjects $S$ of $X$. \end{proof}

\begin{lemma}\label{prop1}
Given $n\geq 0$ and a cardinal $\kappa$ there is a cardinal
$\lambda$ so that if $L$ is an object of injective dimension at most
$n$ and if $S\subset L$ is such that $|S|\leq \kappa$ then there is
an $L'\subset L$ such that
 $S\subset L'$, such that $|L'|\leq \lambda$ and such that both $L'$ and
$L/L'$ have injective dimension at most $n$.
\end{lemma}
\begin{proof} The proof is just the categorical version of [2,
Proposition 2.5]. There the argument uses Corollary 2.3 and Lemma
2.4. To get our version of Corollary 2.3 we only need use the Baer
criterion with respect to a generator $X$ of ${\mathcal{A}}$.
 The categorical version of Lemma 2.4 is just the preceeding lemma.
\end{proof}

\begin{definition}\label{ladefKap}
A class $\Lt$ of objects of ${\mathcal{A}}$ is said to be a
Kaplansky class (see \cite{enochja}) of ${\mathcal{A}}$ if for each
cardinal $\kappa$ there is a cardinal $\lambda$ such that if
$S\subset L$ for some $L\in \Lt$ where $|S|\leq \kappa$ then there
is an $L'\subset L$ with $S\subset L'$ where $|L'|\leq \lambda $ and
where $L'$ and $L/L'$ are both in $\Lt$.
\end{definition}

\begin{corollary}\label{corkap}
If $n\geq 0$ and if $\Lt$ is the class of objects $L$ of
${\mathcal{A}}$ such that $injdim L\leq n$ then $\Lt$ is a Kaplansky
class of ${\mathcal{A}}$.
\end{corollary}
\begin{proof} Immediate from the above. \end{proof}
\medskip

In what follows we will no longer need to refer to a fixed generator
$X$ of $\mathcal{A}$ and so will use the symbol $X$ to stand for an
arbitrary object of $\mathcal{A}$.

\begin{corollary}\label{preenv}
If $n\geq 0$ and if $\Lt$ is the class of objects $L$ of
$\mathcal{A}$ such that $injdim L\leq n$ then $\Lt$ is
preenveloping.
\end{corollary}
\begin{proof} Given the object $X$ of $\mathcal{A}$ we consider
morphisms $X\to L$ where $L\in \Lt$. Using Lemmas
\ref{bob},\ref{cociente} and the Corollaries \ref{setgen} and
\ref{corkap} we see that we can find a cardinal $\lambda$ such that
for any $X\to L$ with $L\in \Lt$ there is a factorization $X\to
L'\subseteq L$ with $L'\in \Lt$ and $|L'|\leq \lambda$. Then using
Corollary \ref{setgen} we see that there is a set $\Lt^*$ of objects
in $\Lt$ such that if $L\in \Lt$ and if $|L|\leq \lambda$ then
$L\cong L^*$ for some $L^*\in \Lt^*$. So noting that for each
$L^*\in \Lt^*$, $\Hom(X,L^*)$ is a set we see that we can find a
family $(\sigma_i)_{i\in I}$ of morphisms $\sigma_i:X\to L_i$ with
$L_i\in \Lt^*$ so that if $L^*\in \Lt^*$ and if $\sigma: X\to L^*$
is a morphism, then $\sigma=\sigma_i$ for some $i\in I$. Then the
morphism $X\to \prod_{i\in I}L_i$ given by the family
$(\sigma_i)_{i\in I}$ is the desired preenvelope.
\end{proof}

\medskip
We remark that for any $X$ of $\mathcal{A}$ we have $X\subseteq E$
for an injective object $E$ of $\mathcal{C}$. Since $E\in \Lt$ we
have the factorization $X\to L\to E$ for any such $\Lt$-preenvelope
$X\to L$. Hence $X\to L$ is necessarily a monomorphism.

\medskip
We will eventually want another property of a class $\Lt$ of
${\mathcal{A}}$. In the next Lemma we will use the notion of
transfinite extensions. For a definition see (\cite[Section
6]{hovey2}).

We recall that since $\mathcal{A}$ is a Grothendieck category (and
so has enough injectives) we can define the functors $\Ext^n(X,Y)$
for all $n\geq 0$.

\begin{definition}
If $X$ is an object of ${\mathcal{A}}$ we say that $projdim\, X\leq
n$ if $$\Ext^i(X,-)=0 \; \mbox{for}\; i\geq n+1.$$ Analogously if an
object $X$ has finite injective dimension at most $n$ we write
$injdim\, Y\leq n$. So then we define $projdim\, X$ and $injdim\, Y$
as usual. We define $FPD({{\mathcal{A}}})$ as the supremum of
$projdim\, X$ over all $X$ such that $projdim\, X<\infty$. We define
$FID({{\mathcal{A}}})$ in a similar manner.
\end{definition}

\begin{definition} A class of objects $\Lt$ of ${{\mathcal{A}}}$ is said to be
closed under transfinite extensions if whenever
$(L_{\alpha})_{\alpha \leq \lambda}$ is a continuous chain of
objects of ${{\mathcal{A}}}$ such that $L_0=0$ and such that
$L_{\alpha +1}/L_{\alpha}\in \Lt$ whenever $\alpha +1\leq \lambda$
we also have $L_{\lambda}$ in $\Lt$.
\end{definition}

\begin{lemma}\label{lemtrans}
Given $n\geq 0$, if $\Lt$ is the class of objects $L$ of
${\mathcal{A}}$ such that $projdim\, L\leq n$ then $\Lt$ is closed
under transfinite extensions.
\end{lemma}

\begin{proof} This result (for modules) is due to Auslander. Our
result follows from Hovey's proof of \cite[Lemma 6.2]{hovey2}. In
that proof he shows that for any object $Y$ of ${\mathcal{A}}$ the
class of $X$ such that $\Ext^1(X,Y)=0$ is closed under transfinite
extensions (this is a categorical version of a theorem of Eklof
\cite{Eklof}). If we let $Y$ range through the class of $n$-th
cosyzygies of objects of ${\mathcal{A}}$ we get the result by using
the fact that $\Ext^{n+1}(X,Z)= \Ext^1(X,Y)$ if $Y$ is such a
cosyzygy for $Z$. \end{proof}

\begin{definition}
We will say that ${\mathcal{A}}$ is a Gorenstein category if the
following hold:

1) For any object $L$ of ${\mathcal{A}}$, $projdim\, L<\infty$ if
and only if $injdim\, L<\infty$.

2) $FPD({{\mathcal{A}}})<\infty$ and $FID({{\mathcal{A}}}) <\infty
$.

3) ${\mathcal{A}}$ has a generator $L$ such that $projdim\,
L<\infty$.
\end{definition}

So now when we say that $({{\mathcal{A}}}, \Lt)$ is a Gorenstein
category we mean that ${\mathcal{A}}$ is such a category and that
$\Lt$ is the class of objects $L$ of ${\mathcal{A}}$ such that
$projdim \, L<\infty$.

If furthermore $FPD({{\mathcal{A}}})\leq n$ and
$FID({{\mathcal{A}}})\leq n$ we will say that $({{\mathcal{A}}},
\Lt)$ has dimension at most $n$ and then define the dimension of
$({{\mathcal{A}}},\Lt)$ to be the least such $n$.

{\bf Remark.} Generators with finite projective dimension in
Grothendieck categories without e\-no\-ugh projectives were also
considered by Hovey in \cite[Section 2]{hovey3}.

\begin{definition}
By a complete projective resolution in ${{\mathcal{A}} }$ we mean a
complex $ {\rm{\bf P}}=(P^n)$ for $n\in \Z$ (so the complex is
infinite in both directions) of projective objects such that
${\rm{\bf P}}$ is exact and such that the complex $\Hom({\rm{\bf
P}},Q)$ is also exact for any projective object $Q$ of
${{\mathcal{A}}}$. If $C$ is the kernel of $P^0\rightarrow P^1$ then
we say that ${\rm{\bf P}}$ is a complete projective resolution of
$C$. An object $C$ is said to be Gorenstein projective if it admits
such a complete projective resolution. Complete injective
resolutions and Gorenstein injective objects are defined dually (see
\cite{EnOv}).
\end{definition}

\begin{definition}
Given a class $\C$ of objects of ${\mathcal{A}}$ then the class of
objects $Y$ of ${\mathcal{A}}$ such that $\Ext^1(C,Y)=0$ for all
$C\in {\C}$ is denoted ${\C}^{\perp}$. Similarly $^{\perp}{\C}$
denotes the class of $X$ such that $\Ext^1(X,C)=0$ for all $C\in
{\C}$.
\end{definition}

\begin{proposition}\label{thickres}
Let $\mathcal{A}$ be a Grothendieck category and let $\C$ ($\J$) be
the class of Gorenstein projective (injective) objects of
$\mathcal{A}$ (there may not be any other than 0). Then $\C^{\perp}$
( $^{\perp}{\J}$) is a thick subcategory of $\mathcal{A}$ containing
all injective and all projective objects of $\mathcal{A}$. So
$\C^{\perp}$ ( $^{\perp}{\J}$) contains all objects of finite
injective dimension and all objects having a finite projective
resolution.
\end{proposition}
\begin{proof} We prove the result for $\C$ and $\C^{\perp}$. A dual
argument gives the result for $\J$ and $^{\perp}{\J}$.

We note that $Y\in \C^{\perp}$ if and only if for every complete
projective resolution ${\rm{\bf P}}$ the complex $\Hom({\rm{\bf P}},
Y)$ is exact. So easily $\C^{\perp}$ is closed under retracts.

If $0\to Y'\to Y\to Y''\to 0$ is an exact sequence in $\mathcal{A}$
and if ${\rm{\bf P}}=(P^n)$ is a complete projective resolution,
then since each $P^n$ is projective we get that $$0\to \Hom({\rm{\bf
P}},Y')\to \Hom({\rm{\bf P}},Y)\to \Hom({\rm{\bf P}},Y'')\to 0$$ is
an exact sequence of complexes. Hence if any two of these complexes
is exact, so is the third. Hence $\C^{\perp}$ is a thick subcategory
of $\mathcal{A}$.

If $E$ is an injective object of $\mathcal{A}$, then since any
complete projective resolution ${\rm{\bf P}}$ is exact,
$\Hom({\rm{\bf P}},E)$ is also exact. So $E\in \C^{\perp}$. If $Q$
is a projective object of $\mathcal{A}$, then by the definition of a
complete projective resolution ${\rm{\bf P}}$ the complex
$\Hom({\rm{\bf P}},Q)$ is exact. Hence we get all the claims about
$\C^{\perp}$. \end{proof}

\begin{definition}
If ${{\mathcal{A}}}$ is Grothendieck with enough projectives define
$Gpd(X)$ (the Gorenstein projective dimension of $X$) in the usual
way. That is, $Gpd(X)=n$ if the first syzygy of $X$ that is
Gorenstein projective is the $n$-th one and $Gpd(X)=\infty$ if there
is no such syzygy. Then define $glGpd({{\mathcal{A}}})$ (the global
Gorenstein projective dimension of ${{\mathcal{A}}}$). Then also
define $Gid(Y)$ and $glGid({{\mathcal{A}}})$ (without assuming
${{\mathcal{A}}}$ has enough projectives).
\end{definition}

\begin{definition}
A pair (${\F}, {\C}$) of classes of objects of ${{\mathcal{A}}}$ is
said to be cotorsion pair if ${\F}^{\perp}={\C}$ and if $^{\perp}{
\C}={\F}$. It is said to be complete if for each $X$ and $Y$ of
${{\mathcal{A}}}$ there exist exact sequences $0\rightarrow
C\rightarrow F\rightarrow X\rightarrow 0$ and $0\rightarrow
Y\rightarrow C'\rightarrow F'\rightarrow 0$ where $F,F'\in {\F}$ and
where $C,C'\in {\C}$.We furthermore say that $({\F}, {\C}) $ is
functorially complete if these sequences can be chosen in a
functorial manner (depending on $X$ and on $Y$) (see Definition 2.3
of \cite{hovey2}).
\end{definition}

We now give our main result.

\begin{theorem}\label{otpri}
If $({{\mathcal{A}}}, \Lt)$ is a Gorenstein category then $(\Lt,
\Lto)$ is a complete and hereditary cotorsion pair on
${\mathcal{A}}$ and $\Lto $ is the class of Gorenstein injective
objects of ${\mathcal{A}}$. If $({{\mathcal{A}}}, \Lt)$ has
dimension at most $n$ then $Gid(Y)\leq n$ for all objects $Y$ of
${\mathcal{A}}$.
\end{theorem}
\begin{proof} To get that $(\Lt, \Lto)$ is a cotorsion pair we only
need argue that $^{\perp}(\Lt^{\perp})=\Lt$. Clearly $\Lt \subset
{^{\perp}(\Lt^{\perp}})$. But since $projdim\, L\leq n$ for all
$L\in \Lt$, $\Lt^{\perp}$ contains all the $n$-th cosyzygies $Y$ of
objects of ${\mathcal{A}}$. If $\Ext^1(L,Y)=0$ for all such $Y$ then
$projdim\, L\leq n$ and so $^{\perp}(\Lt^{\perp})\subset \Lt$ and so
we get $\Lt={^{\perp}}(\Lt^{\perp})$.

As a first step toward arguing that this cotorsion pair is complete,
we want to argue that it is cogenerated by a set, i.e. there is a
set ${\mathcal{S}}$ with ${\mathcal{S}}\subseteq \Lt$ such that
${\mathcal{S}}^{\perp}=\Lt^{\perp}$. But this follows from the fact
that $\Lt$ is closed under transfinite extensions and that it is a
Kaplansky class (see Lemma \ref{lemtrans} and Corollary \ref{corkap}
respectively). For let $\kappa$ be $sup|T|$ where $T$ is the image
of an arbitrary morphism $X\to L$ for $L\in \Lt$ (every such $T\cong
X/Z$ for some $Z\subseteq X$, so we are using Corollary
\ref{epicard}). Now let $\lambda$ be as in Definition \ref{ladefKap}
for this $\kappa$ and $\Lt$. Now let ${\mathcal{S}}$ be a set of
representatives of $S\in\Lt$ such that $|S|\leq \lambda$ (here we
are using Corollary \ref{setgen} with $\kappa$ replaced by
$\lambda$). Then we see that every $L\in \Lt$ can be written as the
union of a continuous chain $(L_{\alpha})_{\alpha\leq \beta}$ (for
some ordinal $\beta$) of subobjects such that $\alpha+1\leq \beta$,
$L_{\alpha+1}/L_{\alpha}$ is isomorphic to an $S\in {\mathcal{S}}$.
Then we appeal to \cite[Lemma 1]{EkTrl} to see that
${\mathcal{S}}^{\perp}=\Lto$.

So now to get functorial completeness we want to appeal to
\cite[Theorem 6.5]{hovey2}. To do so we need to show that $(\Lt,
\Lt^{\perp})$ is small according to \cite[Definition 6.4]{hovey2}.
In this definition, Hovey gives three conditions on a cotorsion pair
in a Grothendieck category that are required for it to be small. In
our situation the cotorsion pair is $(\Lt,\Lto)$. Applied to this
pair (and, again, in our situation) Hovey's conditions are: $i)$
$\Lt$ contains a set of generators of the category, $ii)$
$(\Lt,\Lto)$ is cogenerated by a set ${\mathcal{S}}\subset \Lt$,
$iii)$ for each $L\in {\mathcal{S}}$ there is a given exact sequence
$0\to K\to U\to L\to 0$ such that for any $Y$, $Y\in \Lto$ if and
only if $\Hom(U,Y)\to \Hom(K,Y)\to 0$ is exact for all such exact
sequence.

We have condition $i)$ by our definition of a Gorenstein category.
We have $ii)$ from the above.

We now argue that $(\Lt,\Lto)$ satisfies a slightly weaker version
of $iii)$. We argue that for each $L\in {\mathcal{S}}$ we have a set
of exact sequences $0\to K\to U\to L\to 0$ (one set for each $L\in
{\mathcal{S}}$) so that $Y\in \Lto$ if and only if $\Hom(U,Y)\to
\Hom(K,Y)\to 0$ is exact for all such exact sequences. Our set for a
given $L$ will be a set of representatives of all short exact
sequences $0\to K\to U\to L\to 0$ where $|U|\leq \kappa$ and where
we get the $\kappa$ from Corollary \ref{epicard} when we let the $Y$
of that lemma be our $L$.

So now suppose $G$ is an object such that $\Hom(U,G)\rightarrow
\Hom(K,G) \rightarrow 0$ is exact for all of the exact sequences in
our set.  We want to argue then that $\Ext^1(L,G)=0$ for all $L \in
\Lt$. Because ${\mathcal{S}}$ cogenerates $(\Lt,\Lto)$ it suffices
to argue that $\Ext^1(L,G)=0$ for all $L\in {\mathcal{S}}$. So let
$0\rightarrow G\rightarrow V\rightarrow L\rightarrow 0$ be exact
with $L\in {\mathcal{S}}$. We want to argue that this sequence
splits. But we know that there is a $U\subset V$ such that $|U|\leq
\kappa$ and such that $U\rightarrow L$ is an epimorphism. Then, up
to isomorphism, we can suppose that with $K=G\cap U$, $0\rightarrow
K\rightarrow U\rightarrow L \rightarrow 0$ is one of our sequences.

So consider the commutative diagram

 \DIAGV{50} {0} \n{\ear}  \n{K}
\n{\ear}\n{U}\n{\ear}\n{L}\n{\ear}\n{0} \nn {} \n{}  \n{\sar}
\n{}\n{\sar}\n{}\n{\seql}\n{}\n{}\nn
 {0} \n{\ear}  \n{G}
\n{\ear}\n{V}\n{\ear}\n{L}\n{\ear}\n{0} \diag

Since by hypothesis $K\rightarrow G$ can be extended to
$U\rightarrow G$ we get that the bottom sequence splits.

If $Z$ is a Gorenstein injective object then using a complete
injective resolution of $Z$ it is easy to argue that $\Ext^1(L,Z)=0$
when $projdim\, L \leq  n$. For such a $Z$ is an $n$-th cosyzygy of
some object $W$ of ${\mathcal{A}}$ and
$\Ext^1(L,Z)=\Ext^{n+1}(L,W)=0$. Hence $\Lt^{\perp}$ contains all
the Gorenstein injective objects of ${\mathcal{A}}$.

We now argue that if $G\in \Lt^{\perp}$ then $G$ is Gorenstein
injective. If $0\rightarrow G\rightarrow E^0\rightarrow
E^1\rightarrow \cdots$ is an injective resolution of $G$, then since
$E\in \Lt$ for any injective object $E$ we get that $\Hom(E,-)$
applied  to $0\rightarrow G\rightarrow E^0\rightarrow E^1\rightarrow
\cdots $ gives an exact sequence. So we have the right half of a
complete injective resolution of $G$. We must now show that we can
construct the left half. Since $(\Lt, \Lt^{\perp})$ is complete
there is an exact sequence $0\rightarrow K\rightarrow L \rightarrow
G \rightarrow 0$ with $K\in \Lt^{\perp}$ and with $L\in \Lt$. So
$L\rightarrow G$ is a special $\Lt$-precover of $G$. Let $L\subset
E$ where $E$ is injective. Then $injdim\, E/L<\infty $ and so
$E/L\in \Lt$ and $\Ext^1(E/L,G)=0$. This means that $L\rightarrow G$
can be extended to (necessarily epimorphic) $E\rightarrow G$. Let
$0\rightarrow H\rightarrow E\rightarrow G\rightarrow 0$ be exact.
Then $E\rightarrow G$ is also an $\Lt$-precover. So
$\Hom(\overline{E},E)\rightarrow \Hom(\overline{E},G) \rightarrow 0$
is exact for any injective $\overline{E}$. We now see that
$\Ext^1(\overline{E}, H)=0$. This follows from the exact
$0\rightarrow \Hom(\overline{E}, H)\rightarrow \Hom(\overline{E}, E
)\rightarrow \Hom(\overline{E}, G)\rightarrow \Ext^1(\overline{E},
H)\rightarrow \Ext^1(\overline{E}, E)=0$ and the fact that
$\Hom(\overline{E}, E)\rightarrow \Hom (\overline{E}, G)$ is
surjective. So now replace $G$ with $H$ in the argument above and we
find an analogous $E'\rightarrow H$ with $E'$ injective and such
that this morphism is an $\Lt$-precover of $H$. Since we can
continue this procedure we see that we can construct a complete
injective resolution
$$\cdots \rightarrow E'\rightarrow E\rightarrow E^0\rightarrow
E^1\rightarrow \cdots $$ of $G$. Finally we want to argue that
$Gid(Y)\leq n $ for any object $Y$ of ${\mathcal{A}}$. But if $Z$ is
an $n$-th cosyzygy of $Y$ we have $\Ext^1(L,Z)=0$ for all $L\in
\Lt$. Hence $Z$ is Gorenstein injective.

We now argue that our cotorsion pair $(\Lt,\Lto)$ is hereditary. In
this situation this means that we need argue that if $0\to G'\to
G\to G''\to 0$ is exact with $G',G\in \Lto$, then $G''\in \Lto$. So
we need to argue that $G''$ is Gorenstein injective. We let $L\to G$
be a special $\Lt$-precover of $G$. This means we have an exact
sequence $$0\to H\to L\to G\to 0$$ with $H$ Gorenstein injective. So
we get a commutative diagram \DIAGV{60} {}  \n{}  \n{0}\n{}\n{0}\nn
{}\n{}\n{\sar}\n{}\n{\sar}\nn {}\n{}\n{H}\n{\eeql}\n{H}\nn
{}\n{}\n{\sar}\n{}\n{\sar}\nn{0}
\n{\ear}\n{\overline{G}}\n{\ear}\n{L}\n{\ear}\n{G''}\n{\ear}\n{0}\nn
{}\n{}\n{\sar}\n{}\n{\sar}\n{}\n{\seql}\nn
{0}\n{\ear}\n{G'}\n{\ear}\n{G}\n{\ear}\n{G''}\n{\ear}\n{0}\nn
{}\n{}\n{\sar}\n{}\n{\sar}\nn {}\n{}\n{0}\n{}\n{0}\diag with exact
rows and columns. Then since $\Lto$ is closed under extensions we
get that $L\in \Lto$. So $L\in \Lto\cap \Lt$. Considering an exact
sequence $$0\to L\to E\to M\to 0$$ with $E$ injective we get by
Proposition \ref{thickres} that $M\in \Lt$. So the sequence splits
and we have $L$ injective. We have the exact $0\to \overline{G}\to
L\to G''\to 0$ with $\overline{G}$ Gorenstein injective and where
$\Hom(E,L)\to \Hom(E,G'')\to 0$ is exact if $E$ is injective, since
$\Ext^1(E,\overline{G})=0$ by Proposition \ref{thickres}. So using
the left half of a complete injective resolution of $\overline{G}$
along with the exact $0\to \overline{G}\to L\to G''\to 0$ and an
injective resolution of $G''$ we get a complete injective resolution
of $G''$. \end{proof}

 \vspace{.2in}

\begin{theorem}\label{otro}
If $({{\mathcal{A}}}, \Lt)$ be a Gorenstein category of dimension at
most $n$ having enough projectives. Then for an object $C$ of
$\mathcal{A}$ the following are equivalent:

1) $C$ is an $n$-th syzygy.

2) $C\in \Ltoo$.

3) $C$ is Gorenstein projective.

As a consequence we get that $glGpd({{\mathcal{A}}})\leq n$ and that
$(\Ltoo,\Lt)$ is a complete an hereditary cotorsion pair.
\end{theorem}
\begin{proof} $1)\Rightarrow 2)$ If $C$ is an $n$-th syzygy of $X$
then $\Ext^1(C,L)=\Ext^{n+1}(X,L)$ for any object $L$ of
$\mathcal{A}$. If $L\in \Lt$ then $injdim L\leq n$ so
$\Ext^{n+1}(X,L)=0$. Hence $\Ext^1(C,L)=0$ and so $C\in \Ltoo$.

$2)\Rightarrow 3)$ By Corollary \ref{preenv} we know $C$ has an
$\Lt$-preenvelope $C\to L$. As noted after the proof of that
Corollary, $C\to L$ is a monomorphism. Let $0\to L'\to P\to L\to 0$
be exact where $P$ is projective. Then since $projdim L<\infty$ we
get $projdim L'<\infty$ and so $L'\in \Lt$. Since $\Ext^1(C,L')=0$,
$C\to L$ can be factored $C\to P\to L$. Then $C\to P$ is a
monomorphism and is also and $\Lt$-preenvelope of $C$. So continuing
this procedure we get that we get an exact sequence $$0\to C\to
P^0\to P^1\to P^2\to \cdots$$ where each $P^n$ is projective and
such that if $Q$ is projective then the functor $\Hom(-,Q)$ leaves
the sequence exact.

Now let $$\cdots \to P_2\to P_1\to C\to 0$$ be a projective
resolution of $C$. Since $\Ext^i(C,Q)=0$ for $i\geq 1$ and $Q$
projective we see that $\Hom(-,Q)$ leaves this sequence exact.
Consequently pasting we see that we get a complete projective
resolution of $C$.

$3)\Rightarrow 1)$ is trivial.

The equivalent $1)\Leftrightarrow 3)$ gives that
$glGpd(\mathcal{A})\leq n$.

We know argue that $(\Ltoo,\Lt)$ is a complete cotorsion pair. The
fact that $glGpd(\mathcal{A}\leq n$ gives that for each object $X$
of $\mathcal{A}$ there is an exact sequence $$0\to L\to C\to X\to
0$$ with $C$ Gorenstein projective and $projdim L\leq n-1$ in case
$n\geq 1$ (and with $L=0$ if $n=0$). The argument is essentially the
dual to the proof of \cite[Theorem 11.2.1]{EdO}.

Then using what is called the Salce trick (see \cite[Proposition
7.1.7]{EdO}) we get that for every object $X$ of $\mathcal{A}$ there
is an exact sequence $0\to X\to L\to C\to 0$ with $L\in \Lt$ and
$C\in \Ltoo$. Hence if $X\in (\Ltoo)^{\perp}$ the sequence splits.
So $X$ as a summand of $L\in \Lt$ is also in $\Lt$. So we get that
$(\Ltoo)^{\perp}=\Lt$ and that $(\Ltoo,\Lt)$ is a cotorsion pair.
then completeness follows from the above. Since an exact $0\to L'\to
L\to L''\to 0$ with $L,L''\in \Lt$ gives $L'\in\Lt$ and since
$\mathcal{A}$ has enough projectives we get that $(\Ltoo,\Lt)$ is
hereditary.\end{proof}

\begin{lemma}\label{lem2}
If ${{\mathcal{A}}}$ is Grothendieck with enough projectives then
$$glGpd({{\mathcal{A}}})\leq m\; \Rightarrow
\;FID({{\mathcal{A}}})\leq m$$ and the converse holds if
${{\mathcal{A}}}$ is Gorenstein (so
$glGpd({{\mathcal{A}}})=FID({{\mathcal{A}}})$ in this case). Dually
we have that $glGid({\mathcal{A}})\leq k\Rightarrow FPD({\mathcal{
A}})\leq k$ and the converse holds if $\mathcal{A}$ is Gorenstein
(so $glGid({\mathcal{A}})=FPD({\mathcal{A}})$ in this case).
\end{lemma}
\begin{proof} Let $L$ have finite injective dimension. Now, given
any object $X$ of ${{\mathcal{A}}}$, let $0\rightarrow C\rightarrow
P_{m-1}\rightarrow \cdots \rightarrow P_0\rightarrow X\rightarrow 0$
be a partial projective resolution of $X$. Then $C$ is Gorenstein
projective and so $\Ext^1(C,L)=0$. This gives that
$\Ext^{m+1}(X,L)=0$. Since $X$ was arbitrary, $injdim\, L\leq m$.
Now assume $FID({{\mathcal{A}}})\leq m$ and that ${{\mathcal{A}}}$
is Gorenstein. Let $0\rightarrow C\rightarrow P_{m-1}\rightarrow
\cdots \rightarrow P_0\rightarrow X\rightarrow 0$ be a partial
projective resolution of any $X$. Then if $L\in \Lt$ we have $injdim
\, L <\infty$ so $injdim\, L\leq m$. So $\Ext^{m+1}(X,L)=0$, i.e.,
$\Ext^1(C,L)=0$. So $C\in {^{\bot}\Lt}$. But this means $C$ is
Gorenstein projective by Theorem \ref{otro}. The argument for the
rest of the proof is dual to this argument. \end{proof}

\begin{theorem}
Let ${\mathcal{A}}$ be a Grothendieck category with enough
projectives. Then the following  are equivalent:

1) ${\mathcal{A}}$ is Gorenstein.

2) $glGpd({{\mathcal{A}}})<  \infty$ and
$glGid({{\mathcal{A}}})<\infty$.

Moreover, if (1) (or (2)) holds we have
$$FID({{\mathcal{A}}})=FPD({{\mathcal{A}}})=glGpd({{\mathcal{A}}})=glGid({{\mathcal{A}}}).$$
\end{theorem}
\begin{proof} We have $1)$ implies $2)$ by Theorem \ref{otpri} and
Theorem \ref{otro}.

So now assume $2)$. Let $\Lt$ consist of all object $L$ with
$projdim L<\infty$ and let $\Lt'$ consist of all objects $L'$ with
$injdim L'<\infty$. Now assume $glGpd(\mathcal{A})\leq n$ and that
$glGid(\mathcal{A})\leq n$. If $C$ is an $n$-th syzygy of an object
$X$ of $\mathcal{A}$ then $C$ is Gorenstein projective. So by
Proposition \ref{thickres}, $\Ext^1(C,L)=0$ for any $L\in \Lt$. So
$\Ext^{n+1}(X,L)=0$ for any $X$ and any $L\in \Lt$. So $injdim L\leq
n$ and so $L\in \Lt'$. A dual argument gives that $\Lt'\subseteq
\Lt$ and that if $L'\in \Lt'$ then $projdim L'\leq n$. So we get
that $\Lt'=\Lt$ and so that $({\mathcal{A}},\Lt)$ is a Gorenstein
category.

To get the desired equality we note we have
$FPD({{\mathcal{A}}})=glGid({{\mathcal{A}}})$ and
$FID({{\mathcal{A}}})=glGpd({{\mathcal{A}}})$ by Lemma \ref{lem2}.
So we argue that $glGid({{\mathcal{A}}})=glGpd({{\mathcal{A}}})$. We
use the fact that $({^{\bot}}\Lt,\Lt)$ and $(\Lt,\Lt^{\bot})$ are
complete cotorsion pairs with ${^{\bot}}\Lt={\C}$ the class of
Gorenstein projectives and $\Lt^{\bot}={\J}$ the class of Gorenstein
injectives.

Furthermore we know $\Hom(-,-)$ is right balanced by ${\C}\times
{\J}$ (see \cite[Theorem 1.2.19]{jrlibro}). So we can define
relative derived functors of $\Hom(-,-)$. These are denoted by
$Gext^n(X,Y)$ ($n\geq 0$ for any $X$, $Y$ objects in
${{\mathcal{A}}})$. Once we have this machinery, then we use the
usual argument that the following are equivalent for an abelian
category ${{\mathcal{A}}}$ with enough injectives and projectives
and an $n$ with $0\leq n<\infty$:

1) $projdim\, X\leq n$ for all $X$,

2) $injdim\, Y\leq n$ for all $Y$,

3) $Gext^i(-,-)=0$ for $i\geq n+1$. \end{proof}

\medskip

We get examples of Gorenstein categories with enough projectives by
considering $_{R}\rm{Mod}$ where $R$ is an (Iwanaga) Gorenstein ring
(see 9.1 of \cite{EdO}). These rings are noetherian. But there are
many nontrivial nonnoetherian $R$ such that $_{R}\rm{Mod}$ is
Gorenstein (see \cite{iacobest} and \cite{iacob4}). By a trivial
example we mean one with $l.gldim\, R<\infty$).

\begin{proposition} If ${\mathcal{A}}$ is a Gorenstein category of dimension at most n,
then for every object $Y$ there is an exact sequence
$$0 \rightarrow Y\rightarrow G\rightarrow L\rightarrow 0$$
where $G$ is Gorenstein injective and where $injdim\, L\leq n-1$.
\end{proposition}
\begin{proof} We can mimic the  proof for modules given in
\cite[Theorem 11.2.1]{EdO}. \end{proof} \medskip

We note that if $0\rightarrow Y\rightarrow G\rightarrow L\rightarrow
0$ is as above, then $Y\rightarrow G$ is a special Gorenstein
injective preenvelope of $Y$. Using these preenvelopes we get a
version of relative homological algebra that is called Gorenstein
homological algebra. We see that for each object $Y$ we get a
Gorenstein injective resolution of $Y$. This just means an exact
sequence
$$0\rightarrow Y\rightarrow G^0\rightarrow G^1\rightarrow \cdots$$
where all the $G^n$ are Gorenstein injective and such that $\Hom(-,
G)$ leaves the sequence exact whenever $G$ is Gorenstein injective.
Such a complex is unique up to homotopy  and can be used to give
right derived functors $Gext^i(X,Y)$ of $\Hom$ (these functors were
introduced in \cite{ej1} and were later studied in \cite{Avramov1}
with different notation). There are obvious natural maps
$Gext^i(X,Y)\rightarrow \Ext^i(X,Y)$ for all $i\geq 0$.

The Tate cohomology functors $\widehat{\Ext}^i_{{\mathcal{A}}}
(X,Y)$ (for any $i\in
Z$) are defined as follows.\\
Let $0\rightarrow Y\rightarrow E^0\rightarrow E^1\rightarrow \cdots
\rightarrow
 E^{n-1}\rightarrow G\rightarrow 0$ be a partial injective resolution of $Y$.
 Then $G$ is Gorenstein injective so has a complete injective resolution
which we can take to be
\begin{center}
{${\rm{\bf E}}= \cdots \rightarrow E^{n-1}\rightarrow E^n\rightarrow
E^{n+1}\rightarrow\cdots $}
\end{center}
with $G={\rm Ker}(E^n\rightarrow E^{n+1})$. This complex is unique
up to homotopy and so we define the groups $\widehat{\Ext}^i(X,Y)$
to be the $i$-th cohomology groups of the complex $\Hom(X,{\rm{\bf
E}})$.

If $0\rightarrow Y\rightarrow E'^0\rightarrow E'^1\rightarrow \cdots
$ is an injective resolution of $Y$ then there is a commutative
diagram

\DIAGV{50}
{...}\n{\ear}\n{0}\n{\ear}\n{0}\n{\ear}\n{E'^0}\n{\ear}\n{E'^{1}}
\n{\ear}\n{E'^2}\n{\ear}\n{...}\nn
{}\n{}\n{\sar}\n{}\n{\sar}\n{}\n{\sar}\n{}\n{\sar}\n{}\n{\sar}\n{}\n{}\n{}\nn
{...}\n{\ear}\n{E^{-2}}\n{\ear}\n{E^{-1}}\n{\ear}\n{E^{0}}
\n{\ear}\n{E^1}\n{\ear}\n{E^2}\n{\ear}\n{...} \diagv{0}{0}{0}

The associated map of complexes is unique up to homotopy and gives
rise to natural maps
$$\Ext^i(X,Y)\rightarrow \widehat{\Ext}^ i (X,Y)$$
for $i\geq 0$. In \cite{Avramov1}, Avramov and Martsinkovsky gave a
beautiful result relating these two collections of natural maps. In
their situation they considered finitely generated modules $M$ with
$Gpd(M)<\infty$ over a noetherian ring. And so used a complete
projective resolution as above to get the Tate cohomology functors.
A. Iacob in \cite{iacob1} removed the finitely generated assumption
and also showed how to use the hypothesis $Gid(Y)<\infty$ instead of
the Gorenstein projective dimension hypothesis. Then she showed that
if both $Gpd(X)<\infty$ and $Gid(Y)<\infty$ then the two procedures
give the same groups $\widehat{\Ext}^i(X,Y)$. So due to Iacob's
extension of the Avramov and Martsinkovsky results we get:

\begin{proposition} If ${\mathcal{A}}$ is a Gorenstein category of dimension at most n then
for all objects $X$ and $Y$ of ${\mathcal{A}}$ there exist natural
exact sequences
\begin{center}
{$0\rightarrow Gext^1(X,Y)\rightarrow \Ext^1(X,Y)\rightarrow
\widehat{\Ext}^1 (X,Y)\rightarrow Gext^2(X,Y)\rightarrow \cdots
\rightarrow Gext^n(X,Y)\rightarrow \Ext^n(X,Y)\rightarrow
\widehat{\Ext}^n(X,Y)\rightarrow 0.$}
\end{center}
\end{proposition}

The following is immediate:

\begin{proposition} In the situation above, the following are equivalent:

1) $Gext^1(X,Y)\rightarrow \Ext^1(X,Y)$ is an isomorphism for all
$X$, $Y$.

2) $Gext^i(X,Y)\rightarrow \Ext^i(X,Y)$ is an isomorphism for all
$X$, $Y$ and all $i\geq 1$.

3) $\widehat{\Ext}^i(X,Y)=0$ for all $X$, $Y$ and $i\geq 1$.

4) $\Lt={{\mathcal{A}}}$.

5) ${\J}={\mathcal{I}}$ (${\mathcal{I}}$ is the class of injective
objects).

6) ${\C}={\mathcal{P}}$ (${\mathcal{P}}$ is the class of projective
objects, only in case ${{\mathcal{A}}}$ has enough projectives).

7) $\widehat{\Ext}^i(X,Y)=0$ for all $X$, $Y$ and $i\in \Z$.

8) For a fixed $i$, $\widehat{\Ext}^i(X,Y)=0$ for all $X$, $Y$.
\end{proposition}

\begin{proof} This is just an application of \cite[Theorem
2.2]{hovey2}, (also see \cite[Theorem 8.6]{hovey2}). For the first
claim we use the notation of that paper and let ${\mathcal{
C}}={{\mathcal{A}}}$, ${\mathcal{W}}= \Lt$ and ${\F}=\Lt^{\bot}$.
For the second part let ${\C}={^{\bot}} \Lt$, ${\mathcal{W}}= \Lt$
and ${\F}= {{\mathcal{A}}}$. \end{proof}

\section{The category of quasi-coherent sheaves over ${\bf P^n}(A)$}

\hskip .5cm Let $A$ be a commutative noetherian ring. This section
deals with the category of quasi-coherent sheaves over  ${\bf
P^n}(A)$ and then over closed subschemes $X\subset {\bf P^n}(A)$. We
want to prove that for certain such $X$ this category is Gorenstein.
In this section we will use the fact that $\Qco(X)$ over any scheme
$X$ is equivalent to the category of certain representations of some
quiver $Q$, which may be chosen in various ways, (see \cite[Section
2]{Estr} or \cite[3.1]{hute} for an explanation of this viewpoint).
In case $X\subset  {\bf P^n}(A)$ is a closed subscheme, this quiver
can be taken to have an especially nice form.  First, ${\bf P^n}(A)$
and $\Qco({\bf P^n}(A))$ can be associated with the quiver $Q$ whose
vertices are the subsets
$$v\subseteq \{0,1,2,\cdots,n\},\;\;\;\; v\neq \emptyset ,$$
where there is a unique arrow $v\to w$ when $v\subseteq w$. The
associated ring $\mathfrak{R}$ is such that $\mathfrak{R}(v)$ is the
ring of polynomials with coefficients in $A$ in the variables
$x_i/x_j$, where $0\leq i\leq n$ and $j\in v$. Then, when
$v\subseteq w$, $\mathfrak{R}(v)\to \mathfrak{R}(w)$ is just a
localization (by the multiplicative set generated by the $x_i/x_j$,
$i\in w$, $j\in v$).

{\bf Example}. \rm{If we consider the projective scheme ${\bf
P^1}(A)$, the previous quiver is}$$\{0\}\hookrightarrow
\{0,1\}\hookleftarrow \{1\}$$

\rm{Let us consider the projective scheme ${\bf P^2}(A)$. Then the
corresponding quiver has the form}
\begin{center}
\epsfxsize=5cm \epsffile{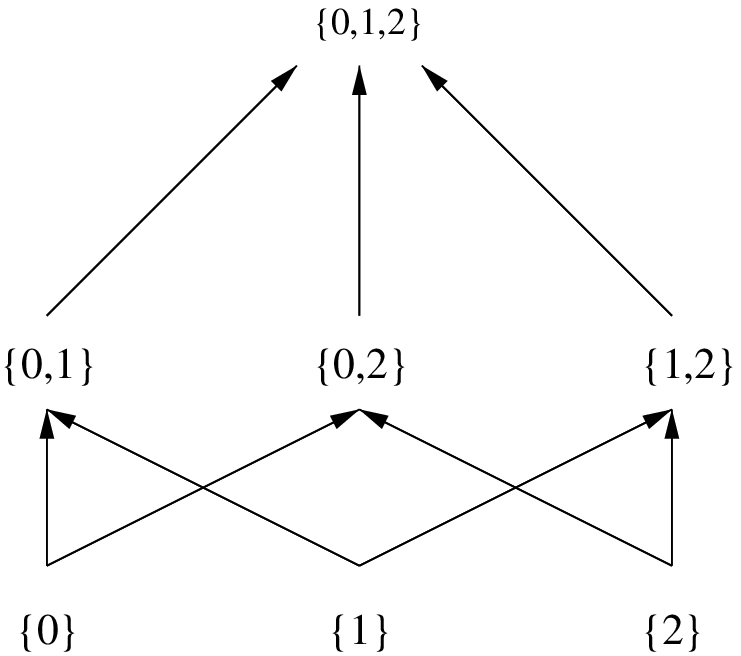}
\end{center}
Notice that, for example, we may delete the arrow from $\{0\}$ to
$\{0,1,2\}$ because this map is the obvious composition
$\{0\}\hookrightarrow \{0,1\}\hookrightarrow \{0,1,2\}$

\medskip\par
A closed subscheme $X\subseteq {\bf P^n}(A)$ is given by a
quasi-coherent sheaf of ideals, i.e. we have an ideal
$\mathfrak{I}(v)\subseteq \mathfrak{R}(v)$ for each $v$ with
$$\mathfrak{R}(w)\otimes_{\mathfrak{R}(v)}\mathfrak{I}(v)\cong \mathfrak{I}(w),$$
when $v\subseteq w$. This just means $\mathfrak{I}(v)\to
\mathfrak{I}(w)$ is the localization of $\mathfrak{I}(v)$ by the
same multiplicative set as above. But then
$$\frac{\mathfrak{R}}{\mathfrak{I}}(v)=\frac{\mathfrak{R}(v)}{\mathfrak{I}(v)}\to
\frac{\mathfrak{R}(w)}{\mathfrak{I}(w)}$$ (when $v\subseteq w$) is
also a localization. So now, to simplify the notation, we will use
$\mathfrak{R}(v)$ in place of
$\frac{\mathfrak{R}(v)}{\mathfrak{I}(v)}$ to give the $\mathfrak{R}$
associated with $X$.

The next result is standard in algebraic geometry. Those who work in
this area will recognize  our $D^v$ as essentially the $i_*$ of
\cite[Proposition II.5.8]{Hartshorne}.
\begin{proposition}\label{eladjunto}
For a given vertex $v$, the functor $$H^v:\Qco(X)\to
{_{\mathfrak{R}(v)} {\mathcal{M}}od}$$ given by $H^v(M)= M(v)$ has a
right adjoint.
\end{proposition}
\begin{proof} We consider $v$ as fixed. Let $N$ be an
$\mathfrak{R}(v)$-module. We construct a quasi-coherent
$\mathfrak{R}$-module $D^v(N)$ as follows: for any $w$ let
$D^v(N)(w)=\mathfrak{R}(v\cup w)\otimes_{\mathfrak{R}(v)}N$ (so
$D^v(N)(w)$ is a localization of $N$). If $w_1\subseteq w_2$ we have
the obvious map
$$D^v(N)(w_1)=\mathfrak{R}(v\cup
w_1)\otimes_{\mathfrak{R}(v)}N\to \mathfrak{R}(v\cup
w_2)\otimes_{\mathfrak{R}(v)}N$$ given by $\mathfrak{R}(v\cup
w_1)\to \mathfrak{R}(v\cup w_2)$. The quasi-coherence of $D^v(N)$
follows from the definition of $D^v(N)$. Given a quasi-coherent
$\mathfrak{R}$-module $M$ we have
$$\Hom_{\Qco(X)}(M,D^v(N))\to \Hom(H^v(M),N).$$ On the
other hand, given $H^v(M)=M(v)\to N$, then for any $w$ we get
$$\mathfrak{R}(v\cup w)\otimes_{\mathfrak{R}(v)}M(v)\to \mathfrak{R}(v\cup
w)\otimes_{\mathfrak{R}(v)}N=D^v(N)(w).$$But we have $M(w)\to
M(v\cup w)$, so composing we get $M(w)\to D^v(N)(w)$ with the
required compatibility. Then it is not hard to check that we have
the required adjoint functor.\end{proof}
\bigskip\par\noindent

{\bf Example.} {\rm We consider $\Qco({\bf P^2}(A))$. Then, if $N$
is an $\mathfrak{R}(\{0,1\})$-module, the quasi-coherent
$\mathfrak{R}$-module $D^{\{0,1\}}(N)$ is given by} \begin{center}
\epsfxsize=4cm \epsffile{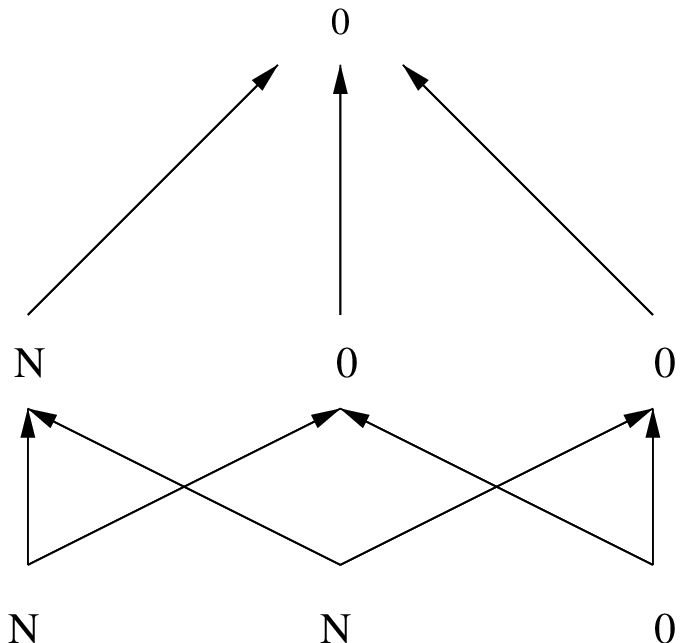}
\end{center}

\medskip\par
The next result is also standard in algebraic geometry (cf. the
comment preceding Proposition \ref{eladjunto}).
\begin{corollary}
With the previous notation, we have:
\begin{itemize}
\item If $N$ is injective then $D^v(N)$ is an injective
$\mathfrak{R}$-module. \item If $N'\to N\to N''$ is an exact
sequence of $\mathfrak{R}(v)$-modules then $D^v(N')\to D^v(N)\to
D^v(N'')$ is also exact.
\end{itemize}
In particular an injective resolution
$$0\to N\to E^0\to E^1\to \cdots$$ gives an injective resolution
$$0\to D^v(N)\to D^v(E^0)\to D^v(E^1)\to \cdots$$
\end{corollary}
\begin{proof} Immediate. \end{proof}

From the previous Corollary we have

\begin{corollary}
Given $N$ and $v$ as above $$\Ext^i_{\Qco(X)}(M,D^v(N))\cong
\Ext^i_{\mathfrak{R}(v)}(H^v(M),N).$$
\end{corollary}
The next Corollary says that if $M$ has finite projective dimension
then ``locally $M$ has finite projective dimension''.
\begin{corollary}\label{converso}
If $projdim\ M<\infty$ then $projdim\ M(v)<\infty$ for every $v$.
\end{corollary}
Now we focus our attention in proving the converse of Corollary
\ref{converso}. To do this we will use a special case of the $N$,
$v$, $D^v(N)$ construction above. Given an $\mathfrak{R}$-module
$M$, $M\neq 0$, choose $v\in Q$ maximal (with respect to
$\subseteq$) so that $M(v)\neq 0$. Then if $v\varsubsetneq w$ we
have $\mathfrak{R}(w)\otimes_{\mathfrak{R}(v)}M(v)=0$. Let $N=M(v)$
and we use this $v$ to construct $D^v(N)$. Now also note that
starting with $M$ and $_{\mathfrak{R}(v)} N$, any $M(v)\to N$ gives
$M\to D^v(N)$. We use this procedure in the special case above where
we chose $v$ maximal such that $M(v)\neq 0$ and let $N=M(v)$. So
letting $M(v)\to N=M(v)$ be $id$ we get a morphism $M\to D^v(N)$ and
so an exact sequence $$0\to K_v\to M\to D^v(N)\to C_v\to
0.$$Moreover, if the previous maximal $v$ is not unique and we
denote by $B$ the set of maximal elements we have the exact sequence
$$0\to K\to M\to \oplus_{v\in B}D^v(N)\to C\to 0.$$
\begin{definition}
By $Supp(M)$ for any quasi-coherent $\mathfrak{R}$-module $M$ we
mean the set of $v\in Q$ such that $M(v)\neq 0$.
\end{definition}
By the construction of $D^v(N)$ from $M$ as above it is easy to see
that $Supp(D^v(N))\subseteq Supp(M)$, so also $Supp(K),\
Supp(C)\subseteq Supp(M)$. But also we see (by the construction)
that $v\notin Supp(K),\ Supp(C)$ since $(M(v)\to D^v(N)(v))=id_N$
for such $v$. So $Supp(K),\ Supp(C)\subsetneq Supp(M)$ when $M\neq
0$.

\bigskip\par\noindent
{\bf Remark.} It follows from the previous definition that
$Supp(M)=\emptyset$ if, and only if, $M=0$. If $|Supp(M)|=1$ then
$Supp(M)$ is of the form $\{i\}$ for some $i$, $0\leq i\leq n$. If
$w\in Supp(M)$ and $w'\subseteq w$ then $w'\in Supp(M)$. If
$|Supp(M)|=1$ and $Supp(M)=\{i\}$ then choosing $v$ as above
(maximal such that $M(v)\neq 0$) we see that $v=\{i\}$ and that with
$N=M(\{i\})$ we have in fact $M=D^v(N)$.
\begin{lemma}\label{reducir}
For any $M$, $projdim\ M<\infty$ if, and only if,
$\Ext^i(M,D^v(N))=0$ for $i>>0$ and for any $N,v$.
\end{lemma}
\begin{proof} The condition is clearly necessary. Now let us assume
the condition. We want to prove $\Ext^i(M,U)=0$ for $i>>0$ and any
quasi-coherent module $U$. If $Supp(U)=\emptyset$ there is nothing
to prove. If $|Supp(U)|=1$, then $U=D^v(N)$ for some $N,v$ and so we
have $\Ext^i(M,U)=0$ for $i>>0$ by hypothesis. So we proceed by
induction on $|Supp(U)|$. But given $U\neq 0$ we construct $$0\to
K\to U\to \oplus_{v\in B}D^v(U(v))\to C\to 0$$ as above (letting $B$
be the set of maximal elements with $U(v)\neq 0$, for all $v\in B$).
Then since $Supp(K),Supp(C)\subsetneq Supp(U)$ we use our induction
hypothesis and easily get $\Ext^i(M,U)=0$ for $i>>0$. \end{proof}
\begin{corollary}\label{projeslocal}
Let $M$ be a quasi-coherent $\mathfrak{R}$-module. Then $projdim\
M<\infty$ if, and only if, $projdim\ M(v)<\infty$ for all $v\in Q$.
Moreover, if $projdim\ M(v)<s$ for all $v\in Q$, then $projdim\
M<s+n$.
\end{corollary}
\begin{proof} The first statement follows by the isomorphism
$$\Ext^i(M,D^v(N))\cong \Ext^i(M(v),N).$$
To see the second one, let $N$ be any  quasi-coherent
$\mathfrak{R}$-module and let $w_N$ be the following non negative
integer: $w_N=max\{ j\;\mid\; Supp(N)$ contains a subset of
cardinality $j\}$.

We consider the exact sequence given in the proof of the  Lemma
\ref{reducir}:
$$0\rightarrow K\rightarrow N\rightarrow \oplus_{v\in B}D^v(N)\rightarrow
C\rightarrow 0,$$ which splits into two short exact sequences:
$$0\rightarrow K\rightarrow N\rightarrow L\rightarrow 0,$$
$$0\rightarrow L\rightarrow \oplus_{v\in B}D_v(N)\rightarrow
C\rightarrow 0.$$ From the second short exact sequence we get the
long exact sequence:
$$\cdots \rightarrow \Ext^{i-1}(M,C)\rightarrow
\Ext^i(M,L)\rightarrow \Ext^i(M,\oplus_{v\in
B}D_v(N))\rightarrow\cdots .$$ Then, we know, by hypothesis, that
$$\Ext^i(M,\oplus_{v\in B}D^v(N))=\oplus_{v\in B}\Ext^i(M(v),N(v))=0$$ for all  $i
> s$.
Hence if $w_N=1$ we have $\Ext^i(M,N)=0$ for all   $i
> s$ because $N$ is a direct sum of $D^v(T)$ for several objects
$T$ and vertices $v$. Then we can prove by induction on $w_N$ that
$\Ext^i(M,N)=0$ for all $i > s+w_N-1$. If $w_N=1$ the result is
proved above. So let $N$ such that $w_N=t$. Then, by the
construction of the exact sequence above, we deduce that $w_C<t$ and
$w_K< t$. Therefore, by induction, $\Ext^{i-1}(M,C)=0$ for all
 $i-1>s+w_C-1$. This implies that
$\Ext^i(M,L)=0$ for all $i> s+w_L-1$ (because $w_L
> w_C$ ). So, from the first short exact sequence, we get the long
exact sequence of $\Ext$: $\cdots\to \Ext^i(M,K)\rightarrow
\Ext^i(M,N)\rightarrow \Ext^i(M,L)\to \cdots $ and, again by
induction applied to $K$, we conclude that $\Ext^i(M,N)=0$ for all
$i> s+w_N-1$ (note that $w_L=w_N$).

If we take $N$ such that $w_N=n+1$, we immediately get that $projdim
M\leq s+n$.
 \end{proof}

\bigskip\par\noindent
{\bf Note.} We also know $injdim\ M<\infty$ if, and only if,
$injdim\ M(v)<\infty$ for all $v$. In fact $injdim\ M=sup_v\ injdim\
M(v)$. As a result of the previous Corollary, the corresponding
statement for $projdim\ M$ is not true.

\bigskip\par\noindent
Now we shall find a family of generators for $\Qco(X)$ with finite
projective dimension. We have the family of $\O(k)$,
$k\in\mathbb{Z}$ for ${\bf P^n}(A)$. These give the family
$\{i^*(\O(k)):\ k\in \mathbb{Z}\}$, where $i:X\hookrightarrow {\bf
P^n}(A)$ (see \cite[pg. 120]{Hartshorne} for notation and
terminology) we will let $\O(k)$ denote $i^*(\O(k))$. Since
$projdim_{\mathfrak{R}(v)}\ \O(k)(v)=0$ for all $k\in \mathbb{Z}$
and all vertex $v$, we get, by Corollary \ref{projeslocal}, that
$projdim\ \O(k)\leq n$ for all $k\in \mathbb{Z}$ (so for an example,
by \cite[Theorem 5.1(c)]{Hartshorne}, $\Ext^n(\O(0), \O(-n-1))\neq
0$ in the $X={\bf P^n}(A)$ case, so we get $projdim\ \O(0)=n$).

It was proved by Serre (see for example \cite{Hartshorne}) that
every coherent sheaf on $X\subseteq {\bf P^n}(A)$ is the quotient of
a finite sum of elements of the family $\{\O(k):\ k\in\mathbb{Z}\}$.
But every quasi-coherent sheaf on $X$ is the filtered union of
coherent subsheafs. So we get that $L=\sqcup_{l\in \Z} \O(k)$ is a
generator for $\Qco(X)$. Furthermore we know that $projdim\ L\leq
n<\infty$. (For a different, and more general, way to get this
family of generators with finite projective dimension see
\cite[Proposition 2.3]{hovey3}).

Now recall that for a Gorenstein ring $B$ (here commutative
noetherian and $injdim\ B<\infty$) we have $$projdim\ L<\infty\
\Leftrightarrow\ injdim\ L<\infty$$ for any $B$-module $L$ (see
\cite{iwa1}). Now suppose that $X\subseteq {\bf P^n}(A)$ is such
that $\mathfrak{R}(v)$ is Gorenstein, for any vertex $v$ ($X$ will
be called, as usual, a locally Gorenstein scheme). Then putting all
of the above together we get.

\begin{theorem}\label{teogor}
If ${{\mathcal{A}}}=\Qco(X)$ for $X\subseteq {\bf P^n}(A)$ a locally
Gorenstein scheme,  then ${\mathcal{A}}$ is a Gorenstein category.
\end{theorem}

\begin{lemma}\label{nose}
Let $M$ be a quasi-coherent $\R$-module, and let $$0\to M\to \E_0\to
\E_1\to \E_2\to \cdots$$ be an injective resolution of $M$. Then
$$0\to M(v)\to H^v(\E_0)\to H^v(\E_1)\to H^v(\E_2)\to \cdots$$ is an
injective resolution of $M(v)$.
\end{lemma}

\begin{proof} It is immediate, because the functor $H^v(-)$ is
exact and preserves injective objects.\end{proof}
\medskip

As a consequence we get that for a quasi-coherent sheaf being
Gorenstein injective is a local property.

\begin{corollary}
Let $X\subseteq {\bf P^n}(A)$ be a locally Gorenstein scheme and $M$
be a quasi-coherent $\mathfrak{R}$-module over $X$. Then $M$ is
Gorenstein injective if and only if $M(v)$ is a Gorenstein injective
$\mathfrak{R}(v)$-module, for all vertex $v$.
\end{corollary}
\begin{proof} We shall use the pair of adjoint functors $(H^v,D^v)$
(with $v\in Q$) obtained in Proposition \ref{eladjunto}.

\medskip\par\noindent
$\Rightarrow )$ Let
$$\cdots\to \E_{-2}\to \E_{-1}\to \E_0\to \E_1\to \E_2\to \cdots$$
be an exact sequence of injective quasi-coherent
$\mathfrak{R}$-modules such that $M=\ker(\E_0\to \E_1)$. Then, for a
fixed vertex $v$, we have an exact sequence of injective
$\mathfrak{R}(v)$-modules
$$\cdots\to H^v(\E_{-2})\to H^v(\E_{-1})\to H^v(\E_0)\to
H^v(\E_1)\to H^v(\E_2)\to \cdots$$ Then if we take an integer
sufficiently large in absolute value, and apply Lemma \ref{nose} and
\cite[Theorem 9.1.11(7)]{EdO} we have that $H^v(\E_{-m})\to
H^v(\E_{-m+1})\to \cdots \to H^v(\E_1)\to H^v(\E_0)\to \cdots$
remains exact when we apply the functor $\Hom(E,-)$, for all $m\geq
0$ and for all injective $\mathfrak{R}(v)$-module $E$ (in fact the
previous is a left ${\mathcal Inj}$-resolution, see page 167 of
\cite{EdO} ). So $M(v)$ is Gorenstein injective.
\medskip\par\noindent
$\Leftarrow )$ Let $M$ be a quasi-coherent $\mathfrak{R}$-module
such that $M(v)$ is Gorenstein injective $\mathfrak{R}(v)$-module.
Since $\Qco(X)$ is a Gorenstein category we may find an exact
sequence $0\rightarrow M\rightarrow G\rightarrow L\rightarrow 0$
with $G$ Gorenstein injective and $L\in {\Lt}.$ Since $M(v)$ and
$G(v)$ are Gorenstein injective (the last by the previous
implication) it follows that $L(v)$ is also Gorenstein injective
($\Lto$ is a coresolving class). Then $L(v)$ is Gorenstein injective
with finite projective dimension, hence an injective module, for all
$v$. Now we take the injective cover of $G$ (which is an epimorphism
with a Gorenstein injective kernel) so we get the exact sequence
$0\rightarrow K\rightarrow E\rightarrow G\rightarrow 0$ with $E$
injective and $K$ Gorenstein injective. Now we make the pull-back
square of $E\rightarrow G$ and $M\rightarrow G$, \DIAGV{60} {}  \n{}
\n{0}\n{}\n{0}\nn {}\n{}\n{\sar}\n{}\n{\sar}\nn
{}\n{}\n{K}\n{\eeql}\n{K}\nn {}\n{}\n{\sar}\n{}\n{\sar}\nn{0}
\n{\ear}\n{U}\n{\ear}\n{E}\n{\ear}\n{L}\n{\ear}\n{0}\nn
{}\n{}\n{\sar}\n{}\n{\sar}\n{}\n{\seql}\nn
{0}\n{\ear}\n{M}\n{\ear}\n{G}\n{\ear}\n{L}\n{\ear}\n{0}\nn
{}\n{}\n{\sar}\n{}\n{\sar}\nn {}\n{}\n{0}\n{}\n{0}\diag Then $U$ is
a quasi-coherent module with finite projective dimension, because it
is part of the exact sequence $0\rightarrow U\rightarrow
E\rightarrow L\rightarrow 0$, and the class $\Lt$ is a resolving
class. Furthermore $U(v)$ is a Gorenstein injective
$\mathfrak{R}(v)$-module, for all $v$, because it is in the middle
of the exact sequence $0\rightarrow K\rightarrow U\rightarrow
M\rightarrow 0$. It follows that $U(v)$ is an injective
$\mathfrak{R}(v)$-module, for all $v$, so $U$ is an injective
quasi-coherent $\mathfrak{R}$-module. So again, since $\Lto$ is
coresolving, we conclude that $M$ is a Gorenstein injective
quasi-coherent $\mathfrak{R}$-module. \end{proof}

\bigskip\par
{\bf Remark.} \rm{It is easy to see that the methods of this section
apply to other schemes. One of the main properties we require of
such a scheme is that the associated quiver $Q$ can be chosen so
that (with the obvious notation) each $\R(v)\rightarrow \R(w)$ is a
localization. This is the case, for example, of toric varieties. We
also point out that the category of quasi-coherent sheaves over
these schemes is equivalent to a quotient category $S-gr/{\mathcal
T}$, for a suitable graded ring $S$ and torsion class ${\mathcal
T}$(see \cite{cox}). This fact could be useful in order to give a
new focus in the topic treated in this paper, taking into account
that the Gorenstein property in graded rings is well-behaved (see
\cite{ja1,ja2}).}

\bigskip\par\noindent
{\bf Acknowledgement.}

The authors are grateful to the referee. He not only suggested
improvements in the paper but also supplied us with a proof that was
lacking.

The final version of this paper was completed during Sergio
Estrada's stay at department of Mathematics of the University of
Kentucky with the support of a MEC/Fulbright grant from the Spanish
Secretar\'{\i}a de Estado de Universidades e Investigaci\'on del
Ministerio de Educaci\'on y Ciencia.

\end{document}